\documentclass[10pt,a4paper]{article}

\usepackage{lscape}

\usepackage{amsmath}
\usepackage{amssymb}
\usepackage{graphicx}
\usepackage{harvard}

\title{
Mass Spectrometry Proteomic Diagnosis: Enacting the Validation
Paradigm. \footnote{The material presented in this paper is
subject to patent application and therefore strictly confidential.
For this reason also,  all explicit references to mass/charge
ranges or
 peaks/peptides identified has had to be removed from the presentation, while the refereeing process continues.
 The authors would however like to include this information in any final published version, once the patent review process has been completed.
 }
 }
\author{
Mertens, Bart J. A.\footnote{Correspondence to: Bart J. A.
Mertens, b.mertens@lumc.nl,
http://www.lumc.nl }, de Noo, M. E.$^+$, Tollenaar, R. A. E. M.$^+$, Deelder, A. M.$^o$ \\
}

\begin{document}
\maketitle \vspace{-0.5cm} \noindent $^*$ Department of Medical
Statistics and Bioinformatics, University of Leiden, P. O. Box
9604, 2300 RC Leiden,
The Netherlands\\
$^+$ Department of Surgery, Leiden University Medical Center, Leiden, The Netherlands.\\
$^o$ Department of Parasitology, Leiden University Medical Center, Leiden, The Netherlands.\\

\begin{abstract}

This paper presents an approach to the evaluation and validation
of mass spectrometry data for construction of an `early warning'
diagnostic procedure. We describe implementation of a designed
experiment and place emphasis on the consistent and correct use of
validation based evaluation - which is a key requirement to
achieve unbiased assessment of the ability of mass spectrometry
data for diagnosis in this setting. Strict adherence to validation
as a scientific principle will however typically imply that the
analyst must make choices. Like all choices in statistical
analysis, validation comes at a cost! We present a detailed and
extensive discussion of the issues involved and propose that much
greater emphasis and requirement for validation would enter
clinical proteomic science.

\end{abstract}

\section{Introduction}
There currently is much interest in application of mass
spectrometry for the construction of new diagnostic proteomic
approaches for the early detection of disease. This is
particularly the case in situations where no reliable diagnostic
tools have yet been developed. 

Generally, in diagnostic research, we may be interested in two key
problems or objectives. The first and foremost of these is to
ascertain whether there is any information in the data (in our
case mass spectra) to allow future cases to be identified with a
high degree of reliability. The second is the identification of
the separating information (in this case proteomic markers)
itself. At this stage, some of our readers may be puzzled, as
naturally, the ability to do the first would imply the presence of
prognostic (marker) information. From a purely statistical
methodological point of view however, the two problems are subtly
different and pose different demands on the analysis. This is
because a proper answer to the first research question requires
rigorous emphasis on fully validated estimation of diagnostic
performance. The second objective on the other hand poses a
feature extraction or definition problem and hence, by its very
nature, may easily cause the analyst to introduce so-called data
`tuning' and optimization steps. This can arise through the need
or desire to employ some form of data cleaning, data preprocessing
steps, use of algorithmic optimization techniques based on the
complete observed data and other similar interventions with the
data. Hence, a conflict may arise and the analyst should carefully
weigh the options and consider the objectives to ensure that the
primary research question can be answered. In particular, if the
first objective is of interest, then special care must be taken to
strictly avoid introducing any optimization steps which could
result in optimistically biased assessment of the diagnostic
ability of the technology, unless some form of validation can  be
introduced to `counteract' and assess the effects. Irrespective of
these considerations, both research questions require carefully
designed experimentation to ensure validity of any study results.

In this paper, we discuss the problem of ascertaining the
viability of using mass spectroscopic analysis of serum samples
for the construction of a diagnostic test for colorectal cancer at
an early stage of the research effort. In other words, we describe
- in essence - a feasibility study. A crucial objective of such
studies is to provide information which allows us to make
decisions as to the continuation of the research effort (which may
involve experiments of much greater cost and complexity in
comparison to the first-stage evaluation). Hence, it is essential
to get a fully validated and unbiased assessment of predictive
error rates and thus the primary research objective is clearly
concerned with the first research question discussed above.

We will now first discuss some additional complicating special
issues which are a concern with mass spectrometry proteomics as
well as a brief discussion of the data, prior to a summary of our
objectives and outline description of the remainder of the paper.

\subsection{Mass spectrometry proteomics, sample size and clinical science}
A key problem in many proteomic studies - but similar problems
arise in many other similar settings ({\it e.g.}: microarray
diagnostics, chemometric discriminant studies) - is the difficulty
to collect a sufficient number of samples. In oncology
applications this will tend to happen, simply because the cancer
of interest may be (relatively) rare. Our example is again a
typical one, as although colorectal cancer is one of the most
common human malignancies, the number of patients any hospital may
expect to encounter on a yearly basis will be limited. On the
other hand, clinicians and biomedical researchers who wish to
explore application of proteomic mass spectrometry technology for
the construction of new diagnostic procedures, will be interested
first to get an indication of whether there is information in the
spectra to allow groups to be separated and what the likely error
rates of misclassification will be. Both these reasons conspire to
cause many proteomic studies to be of small sample size initially.

While small sample size causes problems of precision of model
calibration in the first instance,  it causes special problems in
the proteomics setting because of two complicating reasons. First
is the complexity of the mass spectrometry signal which can
consist of hundreds of possibly overlayed peaks, which are
typically stored on a discretised predefined fine grid of bins. As
a consequence and in combination with small sample size, this
provides the ideal setting for so-called `data-dredging' exercises
which easily generate unsubstantiated claims. The second problem
is that across institutions and research groups we can expect many
such attempts being undertaken to explore application of proteomic
mass spectrometry for diagnosis. For all these reasons, there is
an urgent need for fully validated methods of discriminatory
assessment of proteomic patterns which are sufficiently critical
and conservative and give unbiased estimates of error rate in
small sample size situations.

\subsection{Mass spectrometry data}
The experiment and data discussed and analyses in this paper are
derived from a MALDI-TOF (Matrix Assisted Laser Desorption
Ionisation Time-Of-Flight) mass spectrometer (specifically a
Ultraflex TOF/TOF instrument, Bruker Daltonics, equipped with a
SCOUT ion source which was operated
 in linear mode). The
spectrometer produces a sequence of intensity readings for each
sample on an
ordered set of contiguous bins in the m/z range from UUU 
to VVV
Dalton. Bin sizes (length) of the unprocessed spectra
gradually increase with increasing m/z values, ranging from 0.07
Dalton at the lower end of the mass/charge scale up to 0.24 Dalton
at the upper end of the scale.


 \mbox{}
\\
We now summarize our objectives with this paper, which  are to...
\begin{itemize}
\item[]  ...explain and propose full double cross-validation based
methodology for assessment of diagnostic potential.
\item[] ...explain the limitations which an analyst and researcher
must endure and empose on themselves in order to maintain the
benefits from a fully validatory approach.
\item[] ...propose design suitable for first-stage evaluation of
diagnostic potential and which can also support double
cross-validation.
\item[] ...to describe application to the colon cancer mass
spectral data.
\item[] ...to propose that validation should be given much greater
prominence in wider scientific diagnostic research and that double
cross-validation can play an important role to achieve this in the
first stages of evaluating a new methodology for diagnostic
purposes.
\end{itemize}

We will discuss design first, followed by a description of the
discriminant method and a double cross-validatory approach to
joint estimation and validation of the allocation rule, which
allows for validated error rate evaluation. We pay special
attention to the implied consequences for both the design and
analysis strategy if we are to avoid bias as a consequence of
optimization. As we have explained, this will apply particularly
to preprocessing steps which should be carefully considered to
prevent such problems from affecting the credibility of the first
evaluation of the discriminatory potential of the data. Put
simply: any assessment of such potential must reflect the data,
not the abilities  of the analyst to - perhaps inadvertently -
artificially induce it. In practice, this implies that an analyst
must make choices - perhaps even sacrifices - and we describe the
issues involved. Subsequent to description of the methodological
approach, we discuss application to the colon cancer data and
present a {\it post hoc} exploratory data analysis to
interpretation of the results. While we will focus on our example
to structure the discussion, the issues apply quite generally to
similar problems in proteomics and many other related problems in
bioinformatics, chemometrics, statistical prediction and beyond.
We will assume that the reader has some knowledge of standard
leave-one-out cross-validation.

\section{Design and validation}

\subsection{Design}

A key and characteristic problem of proteomic mass spectrometry
design is the need to cope with the presence of what we may
loosely refer to as so-called `batch effects'. Examples are
plate-to-plate variability,  day-to-day variation and so on. The
presence of such effects is in reality unavoidable. A naive
response to the problem could be to construct designs and
experiments, which avoid or eliminate  the presence of such
effects (assuming that we can anticipate them properly). However,
this will typically imply sample size restrictions on the designs
which are undesirable. A better approach is - rather than try to
eliminate them from the design - to account and accommodate for
these effects in such a way that they do not lead to erroneous or
artificially induced between-group separation. In addition,  this
makes experimentation more realistic also as any `real-life'
application of the methodology would also have to cope with the
presence of such effects.

The problem is a standard example on application of principles of
{\it statistical design} which have been well described in many
classical textbooks and have now been known for almost a century.
The typical first step in applying this knowledge is to first try
and identify - or at least think about - what the potentially
important batch effects might be in advance. In our case plate
variation, but at least also day-to-day variation - such that we
may think of each plate-by-day combination as a batch - commonly
referred to as a `block' in standard statistical terminology. We
may then randomly distribute the available samples from each group
(colon cancer and controls) across the blocks such that
proportions are (as near as) equal within and across blocks for
each group.  For colon cancer, we randomized samples to plates in
such a manner that the distribution of disease stadia was in
approximately equal proportions across plates as well. The
position on the plates of samples allocated to each plate was also
randomized. Each plate was then assigned to a distinct day, which
completes the design. Table~\ref{table0} summarizes the design as
executed on the first week, which provides us with mass spectra on
63 colon cancer patients and 50 healthy controls .

\begin{table}[htb]
\caption{ \normalsize\label{table0} Design as \underline{executed}
on the first week. A replicate of the entire experiment was run on
the subsequent week using plate duplicates. `Stage' refers to the
distribution of cases across the four respective disease stages.
 \normalsize
 }
\begin{center}
\begin{tabular}{l  cccc cc cccc cc cccc c}
\hline\\
Group&\multicolumn{16}{l}{Plates}& Total\\
\\
&1& &&& &&2 & &&& &&3\\
\\

Controls           &17 &&&    & &&17  &&&  &  &&16 &&& & 50\\
Cases              &22 &&&    & &&22  &&&  &  &&19 &&&  & 63\\
\\
\multicolumn{1}{r}{Stage}&\multicolumn{1}{l}{1}&\multicolumn{1}{l}{2}&\multicolumn{1}{l}{3}&\multicolumn{1}{l}{4}&&
&\multicolumn{1}{l}{1}&\multicolumn{1}{l}{2}&\multicolumn{1}{l}{3}&\multicolumn{1}{l}{4}&&
&\multicolumn{1}{l}{1}&\multicolumn{1}{l}{2}&\multicolumn{1}{l}{3}&\multicolumn{1}{l}{4}                        \\
\\
\multicolumn{1}{r}{Cases}&\multicolumn{1}{l}{4}&\multicolumn{1}{l}{10}&\multicolumn{1}{l}{4}&\multicolumn{1}{l}{4}&&
&\multicolumn{1}{l}{4}&\multicolumn{1}{l}{10}&\multicolumn{1}{l}{4}&\multicolumn{1}{l}{4}&&
&\multicolumn{1}{l}{3}&\multicolumn{1}{l}{8}&\multicolumn{1}{l}{4}&\multicolumn{1}{l}{4}                        \\

\\

\\
\hline
\end{tabular}
\end{center}
\end{table}

In our case, it was decided to carry out the experiment in a
single week using three plates only, each of which was assigned to
a consecutive day in the middle of the week - Tuesday to Thursday.
Generally speaking however, it is often wise to try and increase
the number of blocks rather than seek to reduce it when batch
effects are at issue - which is the typical, but actually rather
naive first response of the `uninitiated' in these circumstances.
For example, we could have spread the experiment across several
weeks and use several plates per day with the same design
configuration across weeks, which would allow for a separate
analysis to study and disentangle batch-to-batch variability
(besides the analysis of the main research question). More
importantly, this tends to make the design more robust as the
sample material gets spread across a larger number of blocks which
implies reduced loss of sample material in case problems arise
during the experiment with any specific block (see further on for
an example). An advantage of running the experiment in a single
week is that it may be more easy to maintain the study protocol
across the study, which can be a particular issue if the research
facility carrying out the measurement has to provide for several
customers using the same facilities. Good laboratory management
makes the latter argument less credible however.

The above described design is also referred to as \emph{randomized
block design} in the statistical literature and will ensure that
the batch effects - if they materialize - will not induce an
artificial between-group effect or separation. We refer the reader
to the statistical literature on design of experiments for further
discussion and details of the issues involved, as well as many
other examples of these basic design principles (\cite{Cox},
\cite{Box}, \cite{Neter}, \cite{Fisher1960} among others).

\subsection{Validation}

We can exploit  design to augment cross-validatory analysis. This
is because while sample sizes may be small ({\it i.e.} it is
difficult to get new independent samples), the amount of sample
material available for each sample may be more abundant. This
allows the introduction of so-called replicate samples into the
design.

As the samples are pre-arranged on rectangular plates, a second
`copy' of any plate can be made provided sufficient sample
material is available from each sample. (In our case, sufficient
sample material was available for a second copy only). Thus, we
can duplicate the entire design from the first week and remeasure
the replicate plates through the same design on the second
subsequent week, using new sample material from each sample (but
of course not new samples themselves). With this approach, we thus
generally have available from each $i^{th}$ sample an observation
${\bf x}_i^1=(x_{i1}^1,\ldots,x_{ip}^1)$ of the associated
recorded mass spectrum in the first week, where the vector
elements refer to the measured mass/charge intensities on a
predefined and ordered grid of mass/charges of dimensionality $p$.
In addition, we have for each sample a duplicate measurement ${\bf
x}_i^2=(x_{i1}^2,\ldots,x_{ip}^2)$ obtained from the corresponding
replicate on the corresponding plate measured on the same day one
week later. We may denote the associated class label from each
$i^{th}$ observation as $c(i)$ which takes value in the set of
group indicators $\{1,\ldots,G\}$, where $G$ is the number of
groups. [Note we will drop use of the suffixes 1,2 when the
context makes clear to which week the data relates.]

Unfortunately, the replicate measurements from the third plate are
unavailable due to a technical malfunction which occurred on the
last day of the second week. As a consequence we only have
available the 78 replicates from the first 2 plates in week 2 for
further analysis.

\section{Integrated calibration and validation for classification by double cross-validation}
Given the need to provide fully validated error rates and to keep
the problem of overfitting under control, we will restrict
attention to double cross-validated linear discrimination for
joint calibration and validation \cite{Stone74}. We explain
shrinkage-based estimation and the need for it in linear
discrimination first. Then we explain the double cross-validatory
implementation.

\subsection{Linear classification and shrinkage estimation methodology}
We base classification on Fisher linear discrimination
\cite{Fisher}, which is one of the oldest statistical allocation
methods and certainly the most widely used and successful approach
to statistical classification and pattern recognition to this day.
It has been derived and may be justified based on a variety of
principles of inference, such as maximization of the between-group
separation relative to within-group error in the two-group case
(Fisher's original argument) or the likelihood principle for
normally distributed within-group populations, among others. The
methodology has been amply studied and has been established as an
extremely reliable and robust form of classification and
discriminant analysis - for example, even when assumptions to pool
covariance matrices are not fully met. It is also important to
note that Fisher discrimination does not require an assumption of
within-group normal dispersion. We must refer to the voluminous
literature for further details on the method \cite{Hastie},
\cite{Seber}, \cite{McLachlan}, \cite{Ripley} which is well
established in the applied sciences, such as biology and medicine.
Hastie et al. (2001) contains an up-to-date account of many new
applications which demonstrate the continuing success of the
approach.

Fisher linear discriminant allocation may be defined as assigning
a new observation with feature vector ${\bf x}$ to the group for
which the distance measure
\[
D_g({\bf x})=({\bf x}-\boldsymbol{ \mu}_g) {\bf \Sigma}^{-1} ({\bf
x}-\boldsymbol{ \mu}_g)^T
\]
is minimal, where $g$ denotes the group indicator with $g \in
\{1,\ldots,G\}$, $\boldsymbol{\mu}_g$ the population means and
${\bf \Sigma}$ the population within-group dispersion matrix which
is assumed equal across groups. In practice, the population means
and dispersion matrix will be unknown and hence must be estimated
from the data. In a high-dimensional problem such as in mass
spectrometry proteomics, this leaves us with a difficulty  in
estimating the dispersion matrix as we will typically not be able
to achieve a full rank estimate.

At the risk of some oversimplification of the discussion, there
are basically two ways in which we may remedy the problem so that
the above methodology may again be applied. The first is through
the selection or construction of a set of features which is
reduced in dimensionality, while capturing most of the variability
in the data. In essence, this is the approach which is currently
applied in most of the mass spectrometry proteomics literature (at
the time of writing) through so-called `peak' selection and
construction algorithms, spectral alignment algorithms and other
similar preprocessing procedures, which are applied prior to
subsequent further analysis to construct diagnostic rules, or
other prognostic models. Typical examples are found in papers by
Baggerly, Yasui, Sauve, Morris
(\cite{Baggerly},\cite{Morris},\cite{Sauve},\cite{Yasui1}), among
others. It is important that readers understand that we do not
consider this approach to be fundamentally flawed for mass
spectrometry proteomic data. On the contrary, it is self evident
that mass spectra consist of mixtures of possibly overlaid
intensity peaks corresponding to substances present in the
analyte. Thus, to try to elucidate this structure (first) is in
principle of interest. Problems will however arise when validation
and unbiased assessment of error rates or any other measure of
predictive performance are crucial. This is because preprocessing
may involve
 a high degree of optimization which is itself data-driven. This
applies particularly to problems of peak selection and definition,
and certainly to alignment.


For these reasons and because unbiased error rate estimation is
crucial to our goals, we will pursue the second option, which is
to reduce the preprocessing steps to such an extent that they no
longer introduce any optimization which can no longer be
validated. We then leave the dimensionality of the data intact but
introduce a regularized estimation of the dispersion matrix to
cope with the singularity of the sample dispersion matrix. This
approach can be conservative and may well not provide the most
optimal classifier and error rate, but will at least not be
optimistically biased.

We will explore two distinct forms of regularization, both of
which may be expressed in terms of the canonical decomposition of
the `observed' (or sample) pooled dispersion matrix ${\bf S}={\bf
Q}{\bf \Lambda}{\bf Q}^T $ where ${\bf Q}$ and ${\bf \Lambda}={\rm
diag} (\lambda_1,...,\lambda_{r})$ are the matrices of principal
component weights (or loadings) and variances respectively, with
$\lambda_1>...>\lambda_{r}>0$ respectively (r is the rank of the
pooled covariance matrix). We may now re-estimate the within-group
covariance matrix by only retaining the first $1 \le k\le r$
components only, which gives an estimate
\[
{\bf S}_{(k)}={\bf Q}_{(k)}{\bf
\Lambda}_{(k)}{\bf Q}^T_{(k)},
\]
where ${\bf \Lambda}_{(k)}={\rm diag} (\lambda_1,...,\lambda_{k})$
and ${\bf Q}_{(k)}$ denotes the corresponding reduced matrix of
component loadings. The associated linear discriminant allocation
rule hence assigns observations to the group for which the
smallest sample-based distance estimates
\[
\widehat{D}_g({\bf x})= ({\bf x}-\overline{\bf x}_g) {\bf
S}_{(k)}^{-1} ({\bf x}-\overline{\bf x}_g)^T
\]
are observed, with $\overline{\bf x}_g$ the sample group means for
$g \in \{1,...,G\}$. In the two-group case, this is also
equivalent to least-squares regression analysis using the
Moore-Penrose inverse of  the pooled covariance matrix when $k=r$
(all components kept, also known as shortest least squares
regression), or else ($k<r$) is equivalent to so-called shrunken
least-squares regression (\cite{Ripley} or \cite{Hand97} for more
details). An alternative  is to employ a ridge regularization
\[
{\bf S}_{(\gamma)}={\bf Q}[(1-\gamma){\bf \Lambda}+\gamma {\bf I]}{\bf
Q}^T,
\]
where $0 < \gamma \le 1$ is the ridge regularization or `tuning'
parameter, in which case the sample distance measures are $({\bf
x}-\overline{\bf x}_g) {\bf S}_{(\lambda)}^{-1} ({\bf
x}-\overline{\bf x}_g)^T$.

\subsection{Double cross-validatory estimation and validation}
Application of the above described classification approaches still
requires choice of the tuning parameters $k$ or $\gamma$ involved.
As we are specifically interested in an evaluation of predictive
performance of any diagnostic allocation rule, it becomes crucial
that any optimization - such as the choice of the tuning
parameters - does not take place on the same data used for
validation. On the other hand, predictive tuning is clearly highly
desirable if diagnosis is of interest, so we would not wish to
base the choice of tuning parameters on the full calibration data
itself (and thus effectively drop predictive tuning from the
analysis), but use a truly validatory choice instead. This implies
we either set aside a so-called separate `tuning set' from the
available calibration data prior to validation of predictive
performance itself or appeal to some form of cross-validation.
Good predictive optimization or tuning becomes particulary
important in a high-dimensional setting, such as proteomics, as it
provides a crucial opportunity to safeguard model choice against
overfitting (in other words: over-interpreting the data).
Meanwhile, even if we were able to effectively choose good tuning
parameters, the predictive performance (in our case essentially
the error rates) of any implied allocation rule should again be
validated, which again introduces a need for yet another set-aside
validation set or cross-validation.

We may solve both problems by carrying out a so-called
double-cross-validatory approach, which avoids the need to
introduce separate test (tuning) and validation sets. The method
has been first proposed and investigated by Stone \cite{Stone74}
and integrates predictive optimization and unbiased validated
error rate estimation in a single validatory procedure. While the
principle of the methodology  is sound and well described, this
procedure has until recently not been  applied in practice due to
the considerable computational cost and (algebraic) complexity of
the method. (See  \cite{Mertens03}  for a first full
implementation in the related setting of discriminant allocation
on microarray data.)

As with ordinary leave-one-out cross-validation, double
cross-validation removes each individual (sample) in turn from the
data, after which the discriminant rule is fully recalibrated (and
optimized for prediction) on the leftover data and using the same
procedure in each case. The resulting classification rule is then
applied to the left-out datum to obtain an unbiased allocation for
this sample. This procedure is then repeated across all
individuals and for each person separately, after which
misclassification rates are calculated on the basis of the thus
validated classifications. The double-validatory aspect results
from the fact that the discriminant rule constructed to classify
each left-out datum is optimized through a secondary
cross-validatory evaluation within the first cross-validatory
layer (i.e. full cross-validation again on each `leftover' set
after removal of an observation).  In this manner, we are able to
integrate predictive optimization and predictive unbiased
validation in the same procedure, without loss of data - which is
a crucial requirement to get realistic estimates of error rate
with high-dimensional data.

\section{Application and evaluation}
\subsection{Preprocessing \label{preprocessing}}
Preprocessing is potentially hazardous as it may induce optimistic
bias into the error-rate evaluation. On the other hand,
pre-processing can be beneficial and justified if it removes
variation from the data which does not relate to the group
separation and might obscure an existing group separation. Put
simply, preprocessing is allowed as long as the preprocessing
steps undertaken for any sample are not based on any borrowing of
information or `learning' on the basis of the other samples. In
other words: it must be `within-sample' preprocessing. Note that
this conflicts with some of the preprocessing procedures which
have been applied by other authors in the proteomics literature
(e.g. Coombes, but also Baggerly and Yasui). At this point, we
will describe the preprocessing steps that were employed in our
analysis. We defer a detailed discussion of pre-processing and the
specific potential problems and differences with other authors to
the discussion.

First, we calculated for each sample the average intensity within
each bin across the four mass spectra from the associated spots on
the plate. Then, we aggregated contiguous bins on the m/z scale,
such that the new aggregated bin size spans approximately one
Dalton at the left side of the spectrum and gradually increases to
a width of approximately 3 Dalton at the right hand side. For each
of these new aggregated bins, we calculated for each spectrum the
associated aggregate intensity by summing the intensities across
the bins being aggregated.  Subsequently, spectral baseline was
then removed from each of the thus aggregated spectra separately
using an asymmetric least squares algorithm. 

Suppose ${\bf x}_{b_i}=(x_{b_{i1}},\ldots,x_{b_{ip}})$ denotes the
ordered sequence of baseline corrected m/z intensity values for
the $i^{th}$ sample at this stage of preprocessing. We then
correct the spectrum for the typical intensity and variability
across the spectrum by calculating the standardized values
\[
x_{sb_{ij}}= \frac{x_{b_{ij}}-\text{median}({\bf x}_{b_i})}
{(\text{q}_{0.75}({\bf x}_{b_i})-\text{q}_{0.25}({\bf x}_{b_i}))
},
\]
where  $\text{q}_{0.25}({\bf x}_{b_i})$ and $\text{q}_{0.75}({\bf
x}_{b_i}))$ denote the $25^{th}$ and $75^{th}$ percentiles of the
baseline corrected intensity values for the $i^{th}$ sample. Some
readers will note that these steps bear close resemblance to the
preprocessing procedure proposed by \cite{Satten}, although ours
is a cruder version which does not employ local estimates. The
final preprocessing step is a log-transformation
\[
x_{ij}= \log(x_{sb_{ij}}+\alpha)-\beta
\]
of each spectrum, where $\alpha$ and $\beta$ are two real
constants. We chose $\alpha=100$ and $\beta=4$. The main purpose
of the log-transform is to ensure numerical stability of
calculations.

The above preprocessing steps were applied for each sample and
within each week separately, which thus gives us the observations
${\bf x}_i^1$ and ${\bf x}_i^2$ from the first and second weeks.
It is important to stress that the preprocessing of the data of
any $i^{th}$ sample does not involve use of any information based
on the remaining samples $\{k \mid k \ne i\}$, nor of the
duplicate replicate measured spectrum of the same sample on
another week. This is a vital requirement to ensure the validity
of any cross-validatory evaluation. In addition, there are more
general reasons to avoid any such processing in the
classification/diagnosis context at least, for reasons we will go
in to more detail in the general discussion.

\begin{figure}[tbh]
\centering
\includegraphics[width=15cm]{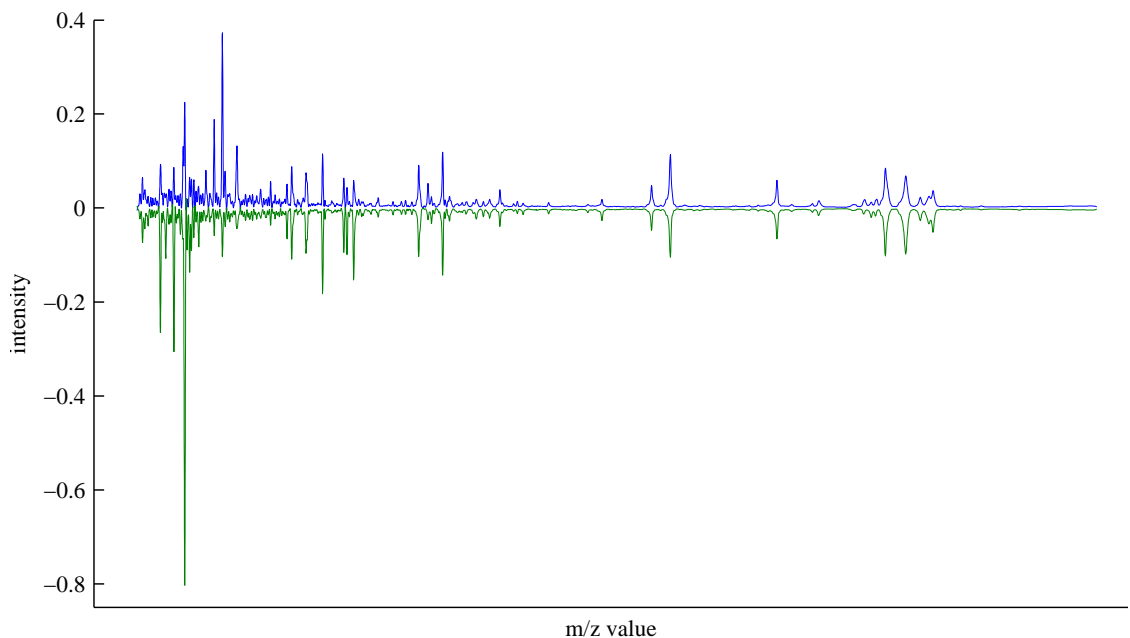}
\caption{Mean spectra for each group separately, after
preprocessing. We plot negative intensity value for the control
group (bottom mean spectrum). } \label{figure1}
\end{figure}

\subsection{Double cross-validatory error rates}
Using the above described preprocessed data, we may carry out a
double- cross-validatory evaluation of predictive performance of
the classifiers we discussed. In the first instance, we restrict
ourselves to the data from the first week. Table~\ref{table1}
 displays the estimated recognition rates and
performance measures from an analysis of the first week data
(leftmost 3 columns). All of the estimates are based on double
cross-validation. We used the average of sensitivity (Se) and
specificity (Sp) as our estimate of the total recognition rate
(T), which implies we assume prior class probabilities to equal
0.5. A threshold of 0.5 was also used to assign observations on
the basis of the a-posteriori class probabilities within the
cross-validatory calculations. B denotes the Brier distance which
we define as
\[
B=\sqrt{  \frac{1}{n} \sum_i{   [1-p(c(i)\mid {\bf x}_i) ]^2}} ,
\]
where $p(c(i)\mid {\bf x}_i)$ is the double cross-validated
predicted a-posteriori class probability for the correct class
$c(i)$ for each $i^{th}$ sample and $n$ is the total sample size.
Likewise, AUC is a double cross-validation estimate of the area
under the empirical ROC curve defined as
\[
AUC={\frac{1}{n_1 n_2}} \sum_{i \in G_1}\sum_{j \in G_2}{[I(p(1
\mid {\bf x}_i)> p(1 \mid {\bf x}_j))+0.5*I(p(1 \mid {\bf
x}_i)=p(1 \mid {\bf x}_j))]},
\]
where $G_1$ and $G_2$ refer to the sample index labels for samples
from the first and second group respectively.

The rightmost three columns of the table refers to a repetition of
this entire double cross-validatory exercise, which replaces each
sample feature vector ${\bf x}_i^1$ with the corresponding
replicate measurement ${\bf x}_i^2$ immediately prior to
classification of that $i^{th}$  sample ({\it i.e.} replacing the
feature vectors with the data from week 2 in the outermost layer
(only!) of the double cross-validatory calculation). Crucially and
importantly, construction of the corresponding discriminant rule
for the classification of each such $i^{th}$ sample in the
internal `calibration' layer of the double cross-validatory
procedure does of course remain based on the data from week 1.
Note that as the replicate data from the third plate are not
available, these results are based on the double cross-validated
predictions for the remaining 78 replicate samples from week 2
only.
\begin{table}[htb]
\caption{\normalsize \label{table1} Double cross-validated
classification results for the colon cancer data. T is the total
recognition rate. Se and Sp are sensitivity and specificity
respectively. B is the Brier distance and AUC is the estimated
area under the ROC curve.} \footnotesize
\begin{center}
\begin{tabular}{llll cc lll}
\hline\\
Method&\multicolumn{3}{c}{First week}& & & \multicolumn{3}{c}{Second week}\\
\cline{2-4} \cline{7-9}
\\
&T (Se,Sp)& B & AUC &&& T (Se,Sp)& B & AUC\\
\\

Moore-Penrose ${\bf S}_{(r)}$ &92.6 (95.2,90.0) & 0.0905  &  97.6                      &&& 94.4 (91.7,97.1)      &0.0885    &97.4 \\
PCA-selection ${\bf S}_{(k)}$ &92.6 (95.2,90.0) & 0.0786  &  97.3                      &&& 88.8 (80.6,97.1)      &0.0935    &96.8 \\
\\
Moore-Penrose Euclidean  ${\bf S}_{(r)}$ ${\bf \Lambda}_{(r)}={\bf I}_{(r)}$   & 89.4 (88.9,90.0) & 0.179  & 96.0           &&& 87.2 (86.1,88.2)      &0.190     &97.0 \\
PCA-selection Euclidisch  ${\bf S}_{(k)}$ ${\bf \Lambda}_{(k)}={\bf I}_{(k)}$  & 88.7 (87.3,90.0) & 0.184  & 96.0           &&& 90.0 (88.9,91.2)      &0.192     &97.0 \\
\\
Ridge    ${\bf S}_{(\gamma)}$                 & 92.0 (95.2,88.0) & 0.0909 & 98.4           &&& 95.8   (91.7,100.0)    & 0.0918   &97.9 \\

\\
\hline
\end{tabular}
\end{center}
\end{table}

At first sight, the Moore Penrose implementation (top line of the
table, both weeks one and two) would seem to be the best
performing and most consistent method. In week 1, Moore-Penrose,
PCA-selection (both using the Mahalanobis distance) and ridge
estimation perform equally well, but there seems to be an increase
in error rate for  week 2 for both the PCA-selection and ridge
implementation. The Euclidean distance based implementations are
worse in the evaluation on the first week, but recognition rates
are consistent  across both weeks when compared to the other
methods.

These results should be interpreted with some caution and require
some explanation. First of all, the `plain' Moore-Penrose is
leave-one-out only as it does not involve choice of shrinkage or
data reduction parameter ($k$ or $\lambda$). The deterioration of
the PCA-selection implementations is partly due to the uncertainty
in estimating the shrinkage terms or choice which is introduced by
the double-cross-validatory estimation.
 For the ridge
implementation, performance is comparable to that from
Moore-Penrose in week 1, which is  not surprising since the chosen
ridge shrinkage parameter $\lambda < 0.0001$ for most
observations. The effects of uncertainty in the determination of
the shrinkage term becomes particularly apparent for PCA-selection
using Mahalanobis distance  (second line in the table) in week 2.

The two Euclidean distance based implementations on the other hand
seem more consistent  across both weeks. The reason is that
component selection is much more stringent for these two
implementations, which selects only the first 2 components for
nearly all observations (with exception of two observations out of
113 for which only the first principal component is retained).
This explains the reduced performance but also the greater
consistency of the classification results. It is precisely because
of this reason that these results (from the Euclidean based
implementations) are more credible and may well turn out to be
more repeatable if the classifier were applied in the future to
data from a new repeat experiment. For comparison, component
selection in the Mahalanobis distance based PCA implementation is
much less stringent and selects ($k=23$ for 53 observations,
$k=28$ for 28 observations and the remainder of the samples uses
even more components). There is thus some evidence of insufficient
shrinkage for this method, and similarly for the ridge
implementation.

\subsection{Investigating bias: a permutation exercise}
We have proposed double cross-validatory integrated estimation and
assessment of statistical diagnostic rules on the basis of the
argument that it should protect against optimistically biased
evaluations. We may check this property by `removing' the class
labels $c(i)$ from the samples $i \in \{1,\ldots,n\}$, randomly
permute and then reassign them to the samples. We then carry out
the double cross-validatory procedure again for any of our
classification methods. Repeating this procedure several times
will give an indication of the biases involved, as the typical
recognition rate  - for example - should equal 50\%  across a
large number of permutations for an unbiased method.

Table~\ref{table2} shows results from such an exercise for the
pca-selection based algorithm across more than 600 such
permutations.
\begin{table}[htb]
\caption{ \normalsize\label{table2} Permutation-based evaluation
of double cross-validatory calculations for linear discrimination
using principal component selection. DBCV refers to the actual
double cross-validatory results (see table~\ref{table1}). $q_{2.5}
$ and $q_{97.5} $ are the 2.5 and 97.5 percentiles. B is the Brier
distance and AUC is the estimated area under the ROC curve.}
\begin{center}
\begin{tabular}{ll lll}
\hline\\
&& \multicolumn{3}{c}{permutation results}\\
\cline{3-5} 
\\
measure &  DBCV &  median  &$\text{q}_{2.5} $ & $\text{q}_{97.5} $ \\
\\

misclassification rate  &7.4    &  50.0  & 36.3    &  72.7  \\
AUC                     &97.3   &  49.4  & 24.8    &  64.2  \\
B                       &0.0786 &  0.324 &  0.200  &  0.446  \\
\\
\hline
\end{tabular}
\end{center}
\end{table}
The results, both for misclassification rate and area under the
(ROC) curve clearly demonstrate the method to be free from bias as
we find median rates and areas of 50\% exactly.

On the other hand, table~\ref{table2} also includes 95\%
confidence intervals for the permutation-based performance
measures. These give an indication of the variability which can be
expected with purely random data and can be compared with the
actually observed double-cross-validation results in our study
(second column of the table). Clearly, the distance between the
validated measures actually observed and even the extreme bounds
of the random permutation confidence intervals is considerable,
which demonstrates the presence of discriminating information in
the mass spectra.

\subsection{Data reduction and post-hoc exploratory analysis}
It is a key feature of our analysis that it places strong emphasis
on estimating fully validated error rates and all steps in the
analysis are geared towards that end. This choice forces a
discipline upon us not to (re-) induce bias via preprocessing,
which specifically precludes any form of either peak `selection'
or `definition'. Therefore, while the analysis is strong on
establishing such unbiased error rates (and other such measures of
predictive performance) and the presence of `predicting' or
`discriminating' information, the analysis is necessarily not
equally outspoken in telling us where this information is to be
found in the spectra. Thus, at this point, we may have convinced
clinicians of the wisdom of randomized block design and (double
cross-) validatory analysis, but still leave an uneasy feeling as
the approach would appear not to be very transparent as to `how'
the classifier assigns observations. We should wish to get an
indication of what the markers are which drive the classification
 - or at least an assurance that the classifier is not
classifying in the noise region of the spectra. Of course, it
would be wrong to label this issue as a deficiency of the
methodological approach. Rather, it is a necessary consequence of
strict adherence to proper protocols for evaluating predictive
ability in a high dimension/small sample size situation.

The dilemma is however not as acute or serious as would appear at
first sight. Indeed, there is nothing stopping us from carrying
out a post-hoc exploratory analysis of both the results and data
(after the validatory calculation). This may include repeating the
analysis after (possibly) some more elaborate preprocessing which
does not have to abide by the limitations of full validation. We
describe two post-hoc exploratory analyses. The first is based on
a very ad hoc algorithmic approach through pre-selection of a
small set of adjacent bins which together account for most of the
variation in the spectra.. The second explores the linear
discriminant weights from a post-hoc fit on the full data.

\subsubsection{Data reduction}
Initialize  $I=\{1,\ldots,p\}$ as the ordered set of bin indices
and $V=\{v_1,\ldots,v_p\}$ the associated set of variances for all
$p$ bins in the preprocessed spectra and across all $n$ samples,
such that $v_j=\sum_i [(x_{ij}-\overline{x}_j)^2]/(n-1)$, where
$\overline{x}_j=\sum x_{ij}/n$ is the sample mean and $j$ is the
bin index number. Calculate the constant
$v_{ref}=\text{q}_{0.95}(V)$ as the 95\% percentile of all $p$ bin
variances. Now initialize the bin selection set $B$ as the set
containing the bin indicator $j$ for which the maximum variance
$v_j$ is observed in the set $V$. Initialize the set of intensity
readings $X_s=\{{\bf x}_{[j]}\mid j\in B\}$ corresponding to the
set $B$, where ${\bf x}_{[j]}=(x_{1j},\ldots,x_{n_j})^T$. We write
${\bf m}=(m_1,\ldots,m_n)^T$ as the set of means
$m_i=\text{mean}(\{x_{ij}\mid j\in B\})$, $i:1,\ldots,n$. Define
$\text{cor}({\bf a},{\bf b})$ to be the coefficient of correlation
between two vectors ${\bf a}$ and ${\bf b}$.

Now run the following algorithm.

\begin{tabbing}

\textit{\{Sta}\=\textit{rt of outer loop\}}\\
    \>\{Sta\=rt of inner loop\}\\

            \>\>Set \=k=1, $I=I-\{j\}$ and $V=V-\{v_{j}\}$\\
            \>\>\>Now  iterate the following procedure until
            termination.\\

            \>\>\>Calculate $\rho_{lower}=\text{cor}({\bf m}_i,{\bf x}_{[j-k]})$ and $\rho_{upper}=\text{cor}({\bf m}_i,{\bf x}_{[j+k]})$\\

            \>\>\>\underline{If} \= $\rho_{lower}>0.9$ and $\rho_{upper}>0.9$ \underline{then}\\

                \>\>\>\>1. Add $j-k$ and $j+k$ to the bin selection set: $B=\{j-k\}\cup{B} \cup \{j+k\}$.\\
                \>\>\>\>2. Update the means $m_i$,$i:1,\ldots,n$.\\
                \>\>\>\>3. \=Remove indices $j-k$ and $j+k$ from the index set $I$, such that $I=I-\{j-k,j+k\}$.\\
                \>\>\>\>\> Similarly update $V=V-\{v_{j-k},v_{j+k}\}$\\
                \>\>\>\>4. set k=k+1\\

            \>\>\>\underline{Else}\\
                \>\>\>\> k=k-1\\
                \>\>\> End iteration.\\

            \> Now select the bin index $j$ for which $v_j=\max(V )$.\\
            \> \underline{If} $v_j>v_{ref}$ \underline{then}\\
            \>\> Update the index set $B=B+\{j\}$ and likewise  $X_s$ and m.\\
            \>\> Go to \{Start of inner loop\}\\

            \> \underline{Else} End algorithm.\\

\end{tabbing}

The algorithm identifies a set of `clusters' of bins. It is
important to note there is no assumption on either shape of the
signal or of monotonicity involved (a single cluster may span
mixture of underlying peaks).  Running this algorithm on the data
from the first week finds the set of indices $B$ that corresponds
to the bins which account for most of the variation in the data.
Applying this to our data results in a subset of 330 bins (in 32
bin clusters  - but it is possible that we visit the same
contiguous region of bins several times). Repeating the entire
double cross-validatory procedure using the principal component
selection shrinkage procedure on this reduced set yields
recognition rates as described in table~\ref{table3}, which are
not inconsistent with those from the full double cross-validatory
evaluation shown in table~\ref{table1}. Hence, we may conclude
that classification is not in the noise region of the spectra.
\begin{table}[htb] \caption{
\normalsize\label{table3} Results from re-running double
cross-validatory calculations after bin-selection for the colon
cancer data (week 1 data only). T is the total recognition rate.
Se and Sp are sensitivity and specificity respectively. B is the
Brier distance and AUC is the estimated area under the ROC curve.}
\normalsize
\begin{center}
\begin{tabular}{llll }
\hline\\

Method&T (Se,Sp)& B & AUC \\
\cline{2-4}
\\
\\

PCA-selection ${\bf S}_{(k)}$ &90.0 (92.1,88.0) & 0.115  &  96.4                     \\
\\
PCA-selection Euclidisch  ${\bf S}_{(k)}$ ${\bf \Lambda}_{(k)}={\bf I}_{(k)}$  & 89.0 (92.1,86.0) & 0.173  & 95.4        \\

\\
\hline
\end{tabular}
\end{center}
\end{table}

\subsubsection{Post-hoc data exploration}
The second aspect which is of interest is a post-hoc exploration
of
 the (linear) discriminant
coefficients ${\boldsymbol\beta}=(\beta_1,\beta_2,...,\beta_p)^T=
{\bf S}_{(k)}^{-1} (\overline{\bf x}_1-\overline{\bf x}_2)^T$ [see
\cite{Seber} or \cite{Hand97}], where $\overline{\bf x}_1$ and
$\overline{\bf x}_2$ are the two sample group means (for cases and
controls). An appropriate and convenient way to summarize and
present the information contained in these coefficients is via the
associated correlations of the measured intensities for each
$j^{th}$ bin with the class indicator, which are easily calculated
as $\rho_j=s_{xj} \beta_j / s_g$, for $j=1,...,p$ where
 $s_{xj}=\sqrt{v_j}$ is the standard deviation at the $j^{th}$
bin and $s_g$ the standard deviation of class indicators.
\begin{figure}[p]
\centering
\includegraphics[width=14cm]{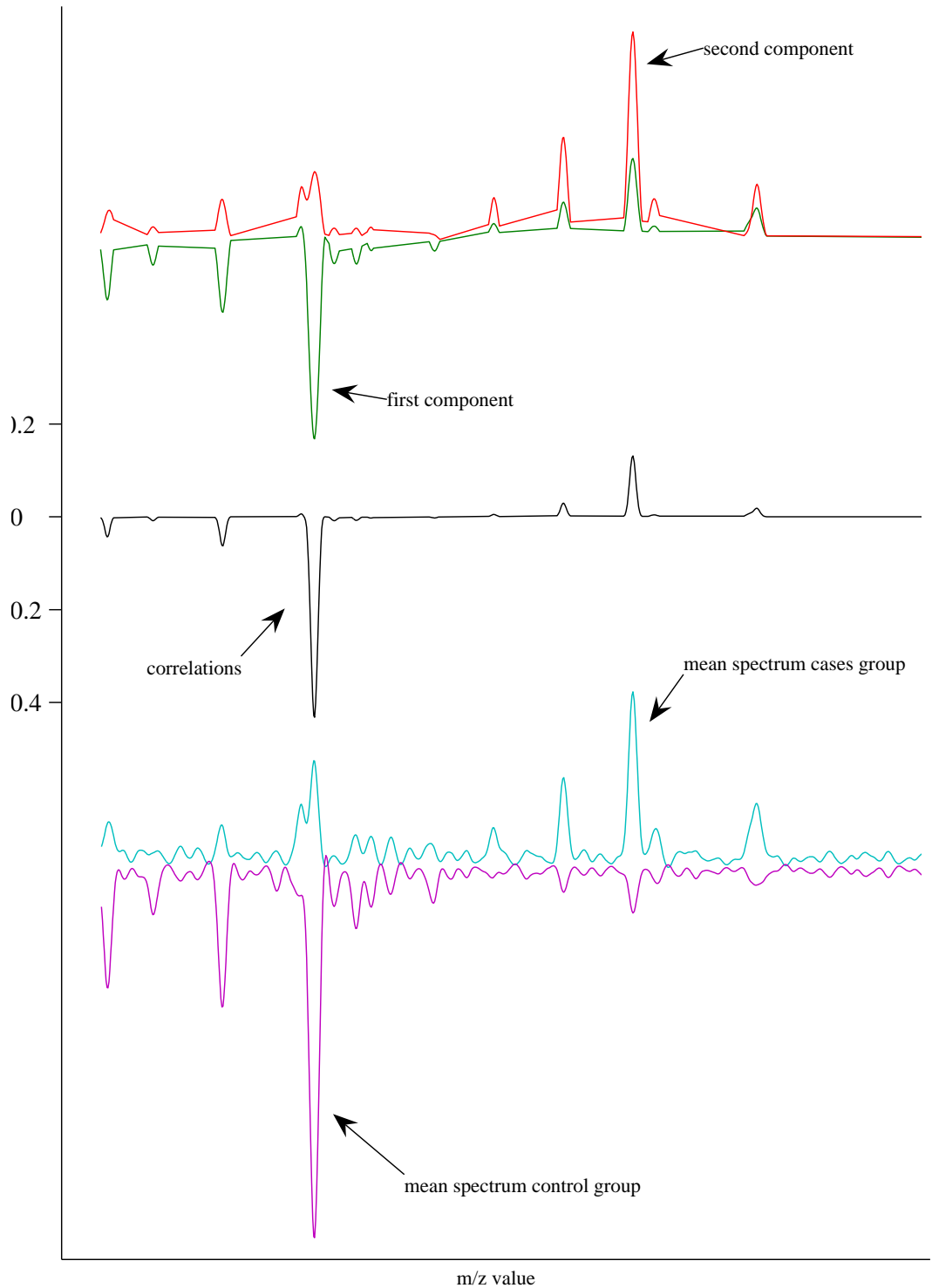}
\caption{Discriminant correlation coefficients $\rho_j=s_{xj}
\beta_j / s_g$ of observed intensity values with the class
indicators in the m/z range from AAA
up to BBB
Dalton. We have plotted the first two principal components above
these correlations for visual comparison and interpretation.
Below the correlations, we plot mean spectra per group ({\it
i.e.}, the vectors $\overline{\bf x}_1$ and $\overline{\bf x}_2$,
as in figure~\ref{figure1}). The y-axis is only relevant to the
correlation coefficient, while we have vertically offset and
rescaled both components and mean spectra to aid visual comparison
across the m/z  range. } \label{figure2}
\end{figure}
We will base this investigation on the linear discriminant fit
using the Euclidean distance on the first two principal components
(use ${\bf S}_{(k)}$, with $k=2$ and ${\bf \Lambda}_{(k)}={\bf
I}_{(k)}$), as   the double validatory assessment of this
classifier clearly identifies the first 2 components as containing
the discriminatory information.

At this point, we can carry out the analysis starting from a
linear discriminant fit based on the full data. Alternatively, we
may equally well base the evaluation on a recomputation of the
linear discriminant fit  on the reduced data described in previous
subsection (in both cases we use the data from the first week).
Figure~\ref{figure2} (middle section) shows a plot of the
correlation coefficients,   subsequent to data reduction
(previously described selection of  330 bins, but of course now
using all 113 samples from the first week). We only show results
within the m/z region between AAA
and BBB
 Dalton, as
 the correlations are effectively zero in the remainder of the m/z
 range. Evidently, this immediately implies that the separating information
 is to be found within the AAA
 to BBB
  m/z range.
We note that the picture shown  is virtually indistinguishable by
eye from that which results from an analysis of the full data (not
shown to save space). The reason for this is that the data
reduction restricts attention to the dominant sources of
variation, which is not very different from what is achieved
through principal component reduction.  Immediately above the
correlation coefficients graph, figure~\ref{figure2} displays the
first two principal components (vertically offset and rescaled to
aid visual interpretation) and again based on the reduced data. In
this case, the distinct bin subsets selected by the previous data
reduction step are clearly visible in the two components, and
display the characteristic `peaks' we would expect to identify.
Disjoint neighboring bin sets are connected with straight lines.
The thus calculated components are a close approximation to those
which would result from an analysis of the full data, as we should
expect (results not shown). As for the correlation coefficients,
any conclusions are therefore identical whether we use the reduced
data or not, although the data reduction step perhaps makes the
component plot easier to `read'.
 At
the bottom of the graph we give the mean spectrum again for each
group separately and from the original data within the m/z range
of interest, as  shown in figure~\ref{figure1} also, along the
complete m/z range.

From this graphical analysis, it is immediately obvious how the
linear discriminant correlation coefficients   identify two major
discriminating contributions, the first of which is centered at
XXX 
Dalton and the second at YYY
Dalton. Furthermore, the correlations have opposite signs at these
locations, which would indicate that the discriminating
information can be summarized through a contrast effect between
corresponding measured intensities in the spectra. An
investigation of the principal components plots above, learns that
the contribution at XXX
Dalton is primarily accounted for by the first component, which
also already contains the contrast with intensities recorded at
YYY
Dalton.  This  contrast is then further amplified by the second
component which identifies a second orthogonal source of variation
relative to the first component, centered predominately at the
already identified peak at YYY
Dalton. Note how each component identifies several other smaller
contributions, which could also be of interest for further
investigation. Comparing these graphs with the within-group mean
spectra, the resemblance with the principal components plots at
the top of the figure are striking and would suggest that the
first component may be primarily explained through variation
within the control group at XXX
Dalton. Likewise, the second
component accounts for a substantial intensity peak  at
Dalton within the colon cancer group.

\begin{figure}[tbh]
\centering
\includegraphics[width=15cm]{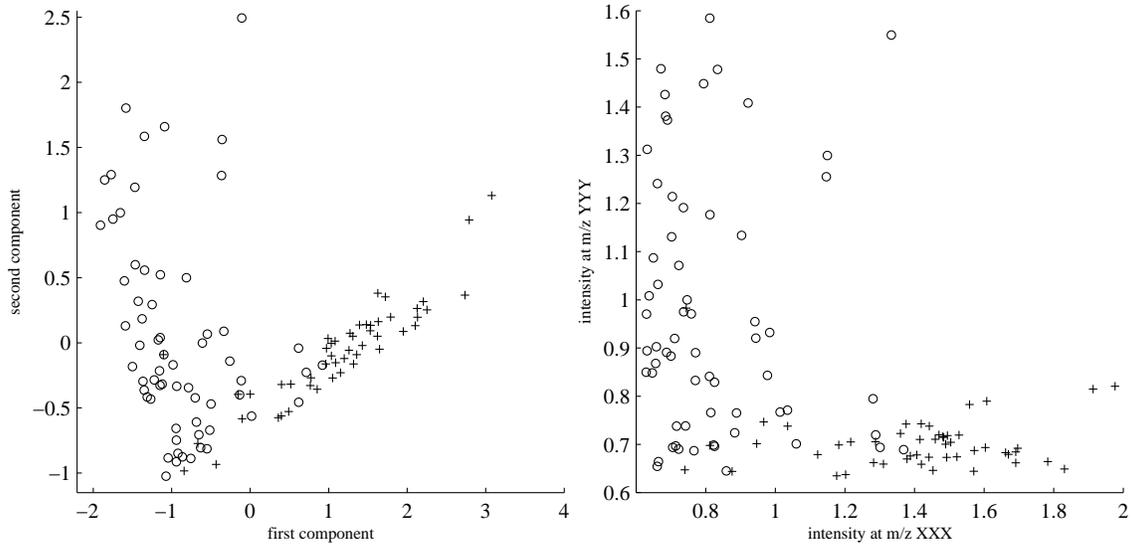}
\caption{Scatter plots distinguishing cases (o) from controls (+).
On the left we plot the second versus the first principal
component. The right plot shows intensity values at YYY
m/z versus those at XXX
m/z. } \label{figure3}
\end{figure}
\begin{figure}[tbh]
\centering
\includegraphics[width=15cm]{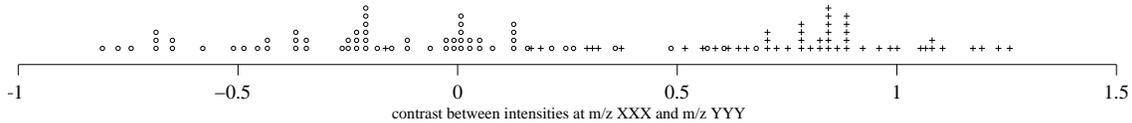}
\caption{Plot of the contrasts (differences) between intensities
at XXX
m/z and
m/z across all observations, using distinct plotting symbols for
each group: cases (o) and controls (+).  } \label{figure4}
\end{figure}
To investigate this further, figure~\ref{figure3} provides scatter
plots of cases and controls versus the first 2 components (left
plot) and between intensities at XXX
and YYY
Dalton respectively (right plot). The resemblance between both
graphs is striking as the right plot can be obtained (virtually)
after clockwise rotation of the left plot. As we can see,
increases in intensity at XXX
Dalton separates controls from cases. Similarly, an increase in
intensity at YYY
Dalton separates cases from controls. The same interpretation
applies to the principal components scatter plot, which confirms
our interpretation of the data in figure~\ref{figure2}.
Figure~\ref{figure4} provides a concise summary graphical
illustration of the results. We calculate the contrast
(difference) for all 113 individuals participating in the study
between the measured intensities at XXX
and YYY
Dalton
and display the differences in a dotplot using distinct plotting
symbols for cases and controls respectively, which demonstrates
the separation between both groups.

For further discussion of the clinical background, study
rationale, setup, execution and interpretation of results from a
substantive clinical perspective, we refer to \cite{Noo} and
subsequent papers from these authors.

\section{Discussion}

\subsection{Current statistical `proteomic' practice}
Most analyses and papers on use of mass spectrometry proteomics
for diagnosis, prognosis and prediction of treatment response are
currently employing a so-called two-stage approach, of which the
first step consists of (extensive) pre-processing procedures and
the second application of some classification or more general
`prediction' algorithm ({\it e.g.} \cite{Baggerly} \cite{Yasui1}
among others). While this approach may at times appear fast on the
way to establishing a `proteomic dogma', there is actually no
agreement as to precisely what pre-processing steps should be
applied and how these may find their place in a coherent analysis
strategy. Thus, while there is much enthusiasm in {\it ad hoc}
algorithm-writing, the methods proposed (and perhaps also results)
may be conflicting. Meanwhile, the dangers of such pre-processing
are little recognized, particularly with respect to the potential
effect on error rate calculation. This may not matter too much if
the objective of the analysis is purely data-exploration. However,
if the objective is to calculate an estimate of prediction error
or the identification of markers that allow us to discriminate
between or prognosticate for \underline{\it future patients}, it
may be a recipe for trouble.

Compare these developments with the analysis presented in
Krzanowski's paper \cite{Krzanowski95} (for example ) on
classification in the related setting of discrimination with near
infrared spectroscopic data. The analysis carried out there is
actually rather cautious and conservative in comparison and does
not require the amount of pre-processing typically presented in
proteomics papers. Nevertheless, Krzanowski is able to demonstrate
presence of an optimistic bias of 10\% in the calculated
recognition rates for one of his classifiers - even though
leave-one-out cross-validation is applied, as a result of the
tuning of only two (!) model parameters based on the full data
(and with a sample size considerable larger than for many current
proteomics projects). In comparison, preprocessing efforts in
proteomics are typically near reckless with extensive
preprocessing steps employed on the full data first, following
which ordinary leave-one-out calculations are frequently presented
as if the first optimization step had never taken place. To make
matters worse, the classifiers employed are often highly complex
and new procedures which must yet pass the test of time and easily
lead to overfitting ({\it e.g.} genetic algorithms
\cite{Baggerly}\cite{Petricoin} ). Good (and standard) statistical
advice however, will be to employ a conservative and inflexible
class of predictors  instead ({\it e.g.} Hand, page 153, or
McLachlan among others), based for example on first and second
moments only.

The extent to which all this will induce optimistically biased
assessments is hard to prejudge, but there can be no doubt that
the current state of affairs will sooner or later induce
false-positive claims on the abilities of mass spectrometry in
specific research projects - if this has not yet materialized.
Long term, the  impact on proteomic science itself may be
disillusionment with proteomics and an averse reaction, as overly
optimistic claims due to improper analysis and study design do not
materialize. A more conservative approach and more stringent
emphasis on full validation in the first stages of evaluation can
help to counteract this, as well as free valuable time and
resources to those projects that are truly worthy of further
attention. In other words, the whole argument for proper
validation and cautious calibration does not restrict itself to
what is `proper' from a purely  statistical point of view - but
also to the economics, viability and credibility of (proteomic)
research and thus to research planning. In this respect, good
statistics leads to good scientific practice.

\subsection{Preprocessing, in particular alignment and feature `definition' }
We want to discuss spectral `alignment' and peak selection in more
detail, as they  pose special problems when applied in any
pre-processing of mass spectra that the analyst should be aware
of, irrespective of the precise methods employed.

Peak `selection' (or discovery as some authors would phrase it) is
really as much peak `definition' as anything else.  From a
practical point of view, the separation of the construction of the
classifier into two phases, the first of which identifies and
aligns the `peaks' may be useful. However, statistical science
 has always regarded the distinction as artificial and more in tune
with the pattern recognition literature \cite{Hand97}(page 51).
Feature extraction and the construction of the classifier both
form part of the optimization of a diagnostic allocation rule. For
that reason, proposals such as by Morris {\it et al.} to base peak
identification on the full data (by means of the average spectrum)
make sense. However as we have already pointed out, it implies use
of the full data and therefore contravenes (cross-)validatory
logic. Another potential problem is that the approach is dependent
on the within-group sample sizes, which may cause difficulties
with unbalanced designs. In principle, cross-validation can be
applied to account for peak-selection as well, but this would
likely be a very costly procedure from a computational point of
view. The other option is to develop new classification methods
which have some form of peak discovery built into the classifier
in such a way that it allows for efficient (cross-)validation (and
for both choice and assessment). If all else fails, use of
separate tuning and validation sets is our only alternative, but
as we have discussed this will often be undesirable in the first
instance due to the sampling costs (see also further on for
further comments on this issue). Even if we were to attempt to
bypass the whole issue by insisting that our only purpose  is
discovery of new `markers' (proteins or peptides), some form of
validation will be of the essence in the construction of the
discriminant models employed as a means to protect against
overfitting.

Alignment is a more tricky problem than peak discovery as it
involves, by definition, reference to other (calibration) samples.
Hence, if cross-validatory analysis is employed within the same
data to either calibrate or validate the classification rule, bias
will be implied by the referencing which has been induced through
aligning previously. At first sight, this might of course be taken
as an argument to avoid cross-validatory evaluation all together
and insist on collecting additional sample material for separate
tuning and test sets. However,  (besides the  increased costs this
will imply) we can not really escape the problems posed by
alignment within classification studies as the issue will return
with a vengeance as soon as we must allocate a new sample: with
respect to which group do we align the new sample for which
diagnosis is required? Perhaps, this could be taken as an argument
for employing so-called `nearest-neighbour type' classifiers (try
alignment of the new sample with - possibly a few samples from -
each group separately). However, these rules are notoriously
unstable in high-dimensional problems. Furthermore, the alignment
problem exists already within the calibration set itself and leads
to strange dilemma's: aligning spectra from one group with those
from another reduces between-group separation (makes them more
alike) and does not seem the right course of action in the
discriminant setting. The opposite choice is equally unappealing
as to align spectra within-group only is sure to induce an
artificial between-group separation.

From a statistical inferential point of view, the whole issue is
therefore more difficult than appears at first sight. While a
proper treatment of the problem on statistical grounds is lacking
(assuming this is possible at all), it might be  advisable not to
align and - in the meantime - place some emphasis on maintaining
good laboratory practice and protocols instead. The only `error'
such choice implies is to make our analyses more conservative
only.

\subsection{Double validatory analysis}
As in our experiment, use of a separate validation set is often
precluded in high dimensional problems, due to sample size
restrictions. Likewise, the need for predictive optimization and
protection against overfitting may be addressed through use of an
additional tuning set, but this greatly increases the burden of
collecting sufficient sample material in practice - often  to the
point of rendering the whole endeavour impossible.

For such reasons, predictive optimization is usually carried out
on the full data instead, which results in optimistically biased
error rate evaluations, particularly with high-dimensional data
such as in mass spectrometry proteomics. The other option is to
reduce the available calibration data prior to optimization so as
to set aside data (perhaps for both a `predictive tuning' as well
as `validation' set) but this approach is not as innocuous as
appears at first sight as it will often reduce the calibration set
beyond what is needed for reasonable calibration. More generally
still, reducing the size of the calibration data changes the
condition of the estimation itself. To put this simply: we are not
only reducing the data by setting-aside aside data from the
calibration set, but also changing the discriminant problem
itself. This is again particularly the case in high-dimensional
cases such as in proteomics where the problem will typically be
ill-conditioned.

The approach we have described in this paper avoids these above
discussed difficulties and dilemmas by reducing the pre-processing
to such an extent that it no longer induces the biases discussed
and then focuses on application of validation for both model
calibration and evaluation through double cross-validation.
Subsequently, a more exploratory analysis can be carried out,
provided we are carefully to interpret results cautiously without
contradicting the primary validated evaluation. We discuss a
number of issues related to application of (double)
cross-validation.

\subsubsection{Full validation}
One potential cause for concern is whether double cross-validation
precludes the need for a completely separate validation set
entirely. Is `double-cross' also `full' validation?

The answer to this question will of course depend on what we mean
by `full validation'. Double cross-validation should give
reasonable protection against overfitting and unbiased estimates
of error rate {\it at the time of  study}. Typically however, the
performance of any decision rule or classifier has a tendency to
`decay' over time. To assess this, subsequent experiments are
needed to verify the estimated error rates. More generally, good
scientific practice requires that we replicate results in a
separate repeat study. This is because cross-validation must
ultimately always remain `within-study' validation and there can
be factors beyond our knowledge which have influenced the study
results. Note however that this applies particularly to the
definition of the case and control group, as the impact of
systematic effects due to measurement can be minimized through use
of randomized block design. Repeat studies may help to detect such
problems. In this sense, double-cross represents the maximum usage
we can make of the data for joint calibration and validation
within a single experiment in order to facilitate decision making
for further follow-up confirmatory experimentation. It is
desirable to have such information available as early as possible,
as subsequent experimentation may be expensive and time consuming.
See also \cite{Ransohoff} for more discussion on these aspects.

\subsubsection{What classifier are we evaluating?}
Two related questions to the previous discussion are `What
classifier does double cross-validation evaluate?' and `How to
assign a new observation?'. Indeed, each observation has its own
classifier in the double cross-validatory evaluation. This seems
to run counter to the intuition that we calibrate a discriminant
rule first and only then evaluate. In that case, the estimated
error rate is taken as a reflection of the diagnostic abilities of
that particular classifier and the allocation of a new sample is
immediate. There is however no logical inconsistency here. Double
cross-validation estimates the error rate we would get `if we were
to apply leave-one-out' on the whole data. Once we know what the
error rate is, we may choose the specific classifier (choice of
$k$ or $\lambda$ in our case) for allocation of future samples (if
required) through application of ordinary leave-one-out on the
whole data (this is in line with the discussion presented by
Mervin Stone \cite{Stone74}).

With double cross-validation, there are however other options to
allow allocation of new samples which have not yet been discussed
in the literature. In our case for example, we may use the mode of
the number of components selected ($k$) across all samples and
then re-estimate the discriminant model with this choice from the
full data. More adventurous still, we could retain each of the $n$
classification rules which are calibrated within the double-cross
procedure and use this ensemble (of classifiers) for allocation of
any future new observation ${\bf x}$. This could be done by
calculating the associated a-posteriori class probabilities
$p_i(g\mid{\bf x})$, for each $i\in\{1,...,n\}$ and
$g\in\{1,...,G\}$, where $p_i$ is obtained from the discriminant
model calibrated in the double-cross procedure when the $i^{th}$
datum has been removed from the data (in the outer shell of the
double-cross procedure). Classification may then be based on the
mean across these $n$ a-posteriori class probabilities for any
$g^{th}$ class. We will not pursue these options further in this
paper.


\subsection{Validation and the future of (statistical) proteomics}
Rigorous emphasis on validation and proper design can help to
establish long-term credibility for proteomic  research and more
general bioinformatics applications. The double-cross approach
with randomized block design described in this paper represents
one contribution towards this goal. Many other steps may however
be taken to enhance the quality of such research studies. One
example is to promote use of `truly' separate validation sets, as
obtained from subsequent separate and additional sampling from the
population of interest and measurement through identical protocols
as applied in the first study. Of course, in practice, this will
only be relevant for those studies which indicate potential from
the first within-study verification of diagnostic ability.

It would be desirable to stimulate greater awareness within the
discipline on the requirements and restrictions which validation
imposes on any analysis.  However, validation and the need for it,
does not only pose a test on the data or any model constructed
from or for it. It does (or should) as much lead to verification
and restriction on the methodological soundness of any proposed
data analysis procedure itself. In this respect we have made
reference to the considerable problems caused by pre-selection and
similarly, so-called alignment methods. The difficulty to
gracefully combine such ad hoc approaches with proper validation
may well indicate a more fundamental problem with the
data-analytic procedures which are used. As a rule of thumb we
should promote forms of analysis which may be validated. This also
implies feasibility of the validation schemes proposed. Of course,
validation is in principle still possible for pre-selection or
construction of markers in principle at least through use of a
separate testing set, but this is not so evident for alignment
procedures. In this paper, we have taken the view that any
methodology should be chosen to err on the conservative side,
whenever such problems arise. Editors of scientific journals can
also contribute much to inspire such conservative attitude by
careful scrutiny of the papers presented for publication. Perhaps
simple check lists could be developed to prevent the major
mistakes from slipping through the net. This may lead to
considerable annoyance in some cases when we face the difficulties
of establishing results in the short term, but may enhance
scientific credibility of the science as a whole in the long run.
Results from the present study show that, with good designed
experimentation, these precautions need not form unsurmountable
obstacles.

%
%
%
%
%
%
%
%
%

\newpage
\section*{Acknowledgements}
The authors would like to give thanks to  Eline Slagboom and Hans
van Houwelingen for considerable support and encouragement in
carrying out this study and to Aliye Ozalp for providing excellent
technical assistance.


\bibliographystyle{agsm}
{}

\end{document}